\documentstyle[11pt]{article} \openup 4pt\oddsidemargin
-.03truein \evensidemargin -.03truein\topmargin -35pt \textwidth
6.66truein\textheight 9.1truein\parskip .09truein \baselineskip
8pt\lineskip 8pt  \def \Z{\hbox{$Z
\hskip -5.2pt Z$}}    \def \C{\hbox{$C\hskip -5pt \vrule height 6pt depth
0pt \hskip 6pt$}} \def\qed{\hfill \hfill \ifhmode\unskip\nobreak\fi
\ifmmode\ifinner \else\hskip5pt\fi\fi \hbox{\hskip5pt\vrule width4pt
height6pt depth1.5pt\hskip 1 pt}} \def\a{\alpha}\def\d
{\delta}\def\D{\Delta} \def\l{\lambda}\def\L
{\Lambda} \def\Si{\Sigma}\def\si{\sigma} \def\VBL{\overline{V}(\L)}\def\sc{\scriptstyle} \def\ssc
{\scriptscriptstyle}\def\dis{\displaystyle} \def\cl{\centerline} \def\nl{\newline}\def\ol{\overline} \def\wh{\widehat}\def\rar{\rightarrow
} \def\lar{\leftarrow} \def\Rla
{\Leftrightarrow}\def\bs{\backslash} \def\rb{\raisebox
} \def\ni{\noindent}\def\hi{\hangindent
}\def\ha{\hangafter} \def\ABS{We classify the finite dimensional
indecomposable $sl(m/n)$-modules with at least a typical or singly
atypical primitive weight. We do this classification not only for
weight modules, but also for generalized weight modules. We obtain that
such a generalized weight module is simply a module obtained by
``linking'' sub-quotient modules of generalized Kac-modules. By
applying our results to $sl(m/1)$, we have in fact completely
classified all finite dimensional $sl(m/1)$-modules.} 
\begin{document}
\cl{\bf CLASSIFICATION OF FINITE DIMENSIONAL MODULES OF SINGLY }\cl{\bf
ATYPICAL TYPE OVER THE LIE SUPERALGEBRAS $sl(m/n)$} \par\ \par\ \par\cl
{Yucai Su\footnote{Partly supported by a grant from Shanghai Jiaotong
University}}\par\ \par\cl{ Department of Applied Mathematics, Shanghai
Jiaotong University, China }\cl{Email: kfimmi@public1.sta.net.cn}\par\
\par\ \par \ni{\bf ABSTRACT.} \ABS\nl{\bf KEYWORDS:} atypical,
primitive, weight diagram, generalized Kac-module, chain. \par\ \par\
\par\ni{\bf I. INTRODUCTION}\par Because finite dimensional
indecomposable modules over Lie superalgebras are not always simple,
the representation theory of Lie superalgebras is more complicated than
that of Lie algebras. Kac in Ref. 1 defined the induced modules $\VBL$
for integral dominant weights $\L$, which are referred to as
Kac-modules by Van der Jeugt {\sl et al} in Ref. 2. Kac divided the
Kac-modules into two categories: {\sl typical} or {\sl atypical}
according as they are simple or not, he also gave a necessary and
sufficient condition for ${\ol V}(\L)$ to be simple. Ref. 2 gave a
character formula for singly atypical simple $sl(m/n)$-modules. Hughes
{\sl et al} in Ref. 3 achieved much progress on the classification of
composition factors of Kac-modules. However, the structure of a
Kac-module is in general still an unsolved problem. More generally, the
problem of classifying finite dimensional indecomposable modules,
posted in Ref. 1, remains open.\par We made a start in Ref. 4, by
giving a complete classification of finite dimensional $sl(2/1)
$-modules. In this paper, we generalize the results to $sl(m/n)$. More
precisely, we classify all finite dimensional $sl(m/n)$-modules with at
least a typical or singly atypical primitive weight. It may be worth
mentioning that although our results here are similar to those in Ref.
4, the proofs are more interesting, more technical, and also more
complicated since Lemma 2.6 in Ref. 4 which was crucial in the proof of
that paper, is no longer valid for general $sl(m/n)$.\par By
introducing the weight diagram, we are able to obtain Theorem 2.9, a
crucial preliminary result in our classification. Then in Sect. III, by
classifying the weight diagram, we obtain our main result of this
section in Theorem 3.8, so that we have a clear picture of a module. In
Sect. IV, by looking deep into generalized weight modules, we can
understand these modules better. Then by a strict and complete proof,
we find out in Theorem 4.9 that such a module is nothing but a module
obtained by ``linking'' some sub-quotient modules of generalized
Kac-modules.\par By applying our results to $sl(m/1)$, we have
efficiently classified all finite dimensional $sl(m/1)$-modules.\par We
would like to point out that it may be possible to use our method to
classify general indecomposable modules as long as we have a better
understanding of the structure of a Kac-module in general.\par\ \par\
\par\ni {\bf II. THE LIE SUPERALGEBRA $sl(m/n)$ AND PRELIMINARY
RESULTS}\par Let $G$ denote the space $sl(m+1/n+1)$ consisting of $(m+n+
2)\times (m+n+2)$ matrices $x=(^{A\ B}_{C\ D})$, where $A\in M_{(m+1)
\times (m+1)}$, $B\in M_{(m+1)\times(n+1)}$, $C\in M_{(n+1)\times(m+1)}
$, and $D\in M_{(n+1)\times(n+1)}$, satisfying the zero supertrace
condition $str(x)=tr(A)-tr(D)=0$. Here, $M_{p\times q}$ denotes the
space of all $p\times q$ complex matrices. Let $G_{\bar0}=\{(^{A\ 0}_{0
\ D})\}$, $G_{\bar1}=\{(^0_C\ ^B_0)\}$, then $G=G_{\bar0}\oplus G_{
\bar1}$ is a $\Z_2(=${}$\Z/2\Z)$ graded space over $\C$ with even part $
G_{\bar0}$ and odd part $G_{\bar1}$. $G$ is a Lie superalgebra with
respect to the bracket relation defined in the above matrix
representation by $[x,y]=xy-(-1)^{ab}yx$, for $x\in G_a$, $y\in G_b$, $
a,b\in\Z_2$. $G_{\bar0}$ is a Lie algebra isomorphic to $sl(m+1)\oplus
\C\oplus sl(n+1)$. Let $G_{+1}=\{(^{0\ B}_{0\ 0})\}$, $G_{-1}=\{(^0_C\ ^
0_0)\}$. Then $G$ has a $\Z_2$-consistent \Z-grading $G=G_{-1}\oplus G_
0\oplus G_{+1}$, $G_{\bar0}=G_0$ and $G_{\bar1}=G_{-1}\oplus G_{+1}
$.\par A Cartan subalgebra $H$ of $G$ has dimension $(m+n+1)$ and
consists of diagonal $(m+n+2)\times(m+n+2)$ matrices of zero
supertrace, with basis $\{h_i=E_{m+i+1,m+i+1}-E_{m+i+2,m+i+2}\,|\,i\ne0,
-m\le i\le n\}\cup\{h_0=E_{m+1,m+1}+E_{m+2,m+2}\}$, where $E_{ij}$ is
the matrix with 1 in $(i,j)$-entry and 0 otherwise. The weight space $H^
*$ is the dual space of $H$ with basis consisting of the simple roots $
\{\a_i\,|\,-m\le i\le n\}$ such that the Chevalley generators are $\{e_
i=E_{m+i+1,m+i+2},f_i=E_{m+i+2,m+i+1}\,|\,-m\le i\le n\}$, and$$\matrix
{\D^+=\{\a_{ij}=\Si_{k=i}^j\a_k\,|\,-m\le i\le j\le n\},\rb{-10pt}{\,}
\hfill\cr \D^+_0=\{\a_{ij}\,|\,-m\le i\le j<0\ \mbox{ or }\ 0<i\le j\le
n\},\ \ \D^+_1=\{\a_{ij}\,|\,-m\le i\le 0\le j\le n\},\hfill\cr}$$ are
respectively the sets of positive roots, positive even roots and
positive odd roots. For $\a_{ij}${}$\in${}$\D^+$, $e_{ij}${}$=${}$E_{m+
i+1,m+j+1}$, $f_{ij}=E_{m+j+1,m+i+1}$ are the generators of the root
spaces $G_{\a_{ij}}$, $G_{-\a_{ij}}$ respectively. Let $(.,.)$ be the
inner products in $H^*$ such that $(\a_i,\a_j)=2$ if $i=j<0$; $=-2$ if $
i=j>0$; $=0$ if $i=j=0$; $=-1$ if $|i-j|=1$ and $i,j\le 0$; $=1$ if $|i-
j|=1$ and $i,j\ge 0$; and $=0$ if $|i-j|\ge2$. Define $\rho_0={1\over2}
\Si_{_{{\sc\a}\in\D^+_0}}\a$, $\rho_1={1\over2}\Si_{_{{\sc\a}\in\D^+_1}}
\a$ and $\rho=\rho_0-\rho_1$. We give a well order in $H^*$: for $\l,
\mu\in H^*$, $\l>\mu\ \Rla\ \l-\mu=\Si_{i=-m}^n a_i\a_i$ such that for
the first $a_i\ne0$, we have $a_i>0$.\par Now let $\L$ be an integral
dominant weight over $G$, i.e., $\L(h_i)\in\Z_+$ if $i\ne0$. Let $V^0(
\L)$ be the simple highest weight $G_0$-modules with the highest weight
$\L$. As in Ref. 4, we first give some definitions and preliminary
results. \par\ni{\sl Definition 2.1.} Extend $V^0(\L)$ to a $G_0\oplus
G_1$-module by requiring $G_1 V^0(\L)=0$, and define {\sl Kac-module}
to be the induced module ${\ol V}(\L)=\mbox{Ind}^G_{G_0\oplus G_1}V^0(
\L)=U(G)\otimes_{U(G_0\oplus G_1)}V^0(\L)\cong U(G_{-1})\otimes V^0(\L)
$. Similarly, one can define the anti-Kac-module ${\ol V}_*(-\L)$ with
the lowest weight $-\L$ in the obvious way, starting from the lowest
weight $G_0$-module $V^0_*(-\L)$ with the lowest weight $-\L$.\qed\par
It follows that any highest weight module with highest weight $\L$ is a
quotient of $\VBL$. We will denote $V(\L)$ the simple highest weight
module with highest weight $\L$. \par\ni{\sl Definition 2.2.} Let $\L$
be integral dominant. $\L$, $\VBL$, $V(\L)$ are called {\sl typical} if
there does not exist an odd root $\a_{ij}\in\D^+_1$ such that $(\L+\rho,
\a_{ij})=\Si_{k=i}^0\L(h_k)-\Si_{k=1}^j\L(h_k)-i-j =0$, otherwise they
are called {\sl atypical}. They are called {\sl singly atypical} if
there exists exactly one such $\a_{ij}$ (in this case, $\a_{ij}$ is
called an {\sl atypical root} of $\L$), otherwise they are called {\sl
multiply} atypical (see Refs. 1, 2 and 3). They are called
anti-typical, singly (multiply) anti-atypical if $-\L$ is typical,
singly (multiply) atypical respectively.\qed\par\ni{\sl Definition
2.3.} A weight vector $v_\l\ne0$ in a module $V$ is called {\sl
primitive} (Ref. 5, \S9.3) if there exists a submodule $U$ such that $v_
\l\not\in U$ but $G_+v_\l\subset U$. If $U=0$, i.e., $G_+v_\l=0$, then $
v_\l$ is called {\sl strongly primitive}. Correspondingly, $\l$ is
called primitive, or strongly primitive. Denote by $P_V$ the set of
primitive weights of $V$. Similarly, one can define anti-primitive
vectors and weights by replacing $G_+$ by $G_-$ in the definition.\qed
\par In this paper, (anti-) primitive vectors are restricted to be
those which are $G_0$-strongly (anti-) primitive and which generate
indecomposable submodules; different primitive vectors always means
they generate different submodules.\par In this paper, we will only
consider those finite dimensional indecomposable modules with a
condition that there exists at least a primitive weight of typical or
singly atypical type (note that such condition does not necessarily
imply that all primitive weights are typical or singly atypical, we
will prove this implication in Theorem 2.9). As $H$ does not always act
diagonally on a $G$-module, in Sects. II and III, we will first
consider weight modules $V$, i.e., $V$ admits a weight space
decomposition: $V=\oplus_{\l\in H^*}V_\l$, where $V_\l=\{v\in V\,|\,hv=
\l(h)v$, for $h\in H\}$. (In Ref. 4, $V$ is called a module with
diagonal Cartan subalgebra.) Then in Sect. IV, we will extend our
results to generalized weight modules.\par\ni {\sl Lemma 2.4.} (1) $
\VBL$ is simple $\ \ \Rla\ \ \L$ is typical.\nl (2) Suppose $\VBL$ is
not simple, then (i) if $\L$ is singly atypical, $\VBL$ has two
composition factors: one is $V(\L)$, the other we will denote by $V(\L^-
)$; (ii) if $\L$ is multiply atypical, all primitive weights of $\VBL$
are multiply atypical.\par \ni{\sl Proof.} (1) See Ref. 1. (2)(i) See
Refs. 2 and 3. (ii) (See Ref. 3) If $\L$ is multiply atypical, then the
bottom composition factor of $\VBL$ is the only simple submodule of $
\VBL$, and its highest weight is multiply atypical. By (i), a highest
weight module of singly atypical type (which has at most two
composition factors) does not contain such a submodule. Thus a simple
highest weight module of singly atypical type cannot be a composition
factor of $\VBL$.\qed\par\ni{\sl Remark 2.5.} If $\L$ is singly
atypical, by introducing in Ref. 3 the {\sl atypicality matrix} $A(\L)
$, an $(m+1)\times(n+1)$ matrix with $(m+i+1,j+1)$-entry being $A(\L)_
{ij}=(\L+\rho,\a_{ij})$, and the southwest chain of $A(\L)$, one can
obtain $\L^-$ by subtracting from $\L$ those $\a_{ij}$ sitting on the
chain. This can be simply done as follows: choose $\L_0=\L$, $\L_1=\L_0-
\a_{ij}$, where $\a_{ij}$ is the atypical root of $\L$, i.e., $A(\L_0)_
{ij}=0$; suppose we have chosen $\L_k=\L_{k-1}-\a_{ij}$, whose
atypicality matrix $A(\L_k)$ is obtained from $A(\L_{k-1})$ by
subtracting 1 from $(m+i+1)$-th row and adding 1 to $(j+1)$-th column,
then $\L_{k+1}=\L_k-\a_{i'j'}$, where $(i',j')\ne(i,j)$ is the another
entry satisfying $A(\L_k)_{i'j'}=0$; continue this procedure until
there is no more such $(i'j')$, then $\L_k$ is dominant and $\L^-=\L_k
$.\qed\par Suppose $\L$ is singly atypical. We will use the following
notations through the paper.\par\ni{\sl Definition 2.6.} (1) Let $\L^-$
denote the weight defined by Lemma 2.4. Suppose $v_\L$, $v_{\L^-}\in
\VBL$ are the primitive vectors, then there exists $g\in G_{-1}U(G^-)$
such that $g v_\L=v_{\L^-}\,$. Fix such a $g$ and denote it by $g^-_\L
$.\nl (2) Let $\L_{low}$ denote the lowest weight in $V(\L)$.\nl (3)
Let $\L^+$ denote the highest weight of the anti-Kac-module ${\ol V}_*(
\L_{low})$. (It follows that $\L_{low}$ is the lowest weight of the
Kac-module ${\ol V}(\L^+)$.)\nl (4) In anti-Kac-module ${\ol V}_*(\L_
{low})$, where $\L_{low}$ is a singly anti-atypical lowest weight,
similar to Lemma 2.4(2.i), there are 2 primitive vectors $v_\L$, $v_{\L^
+}$ in ${\ol V}_*(\L_{low})$. By decomposing the enveloping algebra $U(
G)=U(G^-)U(G_{+1})U(G^+_0)$, and by noting that $v_{\L^+}$ is the
highest weight vector in ${\ol V}_*(\L)$, there exists $g\in G_{+1}U(G_
{+1})$ such that $g v_\L=v_{\L^+}\,$. Fix such a $g$ and denote it by $
g^+_\L$.\qed\par By the above definition, we see from Refs. 2 and 3
that $\L=(\L^+)^-=(\L^-)^+$ and we can compute $\L^+$ as in Remark 2.5
by defining the northeast chain of $A(\L)$ and adding to $\L$ those $\a_
{ij}$ sitting on the chain, so that $\L^+$ is the last $\L_k=\L_{k-1}+
\a_{i'j'}\,$. Thus, for a singly atypical weight $\L$, we can define
inductively $\{\L^{(i)}\,|\,i\in\Z\}$ by $\L^{(0)}=\L$, $\L^{(-i)}=(\L^
{(-i+1)})^-$, $\L^{(i)}=(\L^{(i-1)})^+$, $i>0$. We have $(\L^{(k)})^+=
\L^{(k+1)}$ for $k\in\Z$. For each $\L$, denote $\phi_\L=\{\L^{(i)}\,|
\,i\in\Z\}$, and for $i,j\in\Z$, $i\le j$, denote $\phi_\L^{(ij)}=\{\L^
{(k)}\,|\,i\le k\le j\}$.\par In the following, we will see a module is
uniquely determined, up to an isomorphism, by the relationship between
its primitive vectors. We define a diagram to express the structure of
a module $V$, where two primitive vectors $v_\l,v_\mu$ are linked by a
line with an arrow: $$\mbox{(i) }v_\l\rar v_\mu\ \Rla\ v_\mu\in G_{-1}U(
G^-)v_\l\ \mbox{ and (ii) }v_\l\lar v_\mu\ \Rla\ v_\l\in G_{+1}U(G_{+1})
v_\mu.\eqno(2.1)$$ It follows that a primitive vector $v_\l$ can be
linked by 4 ways: $\rar v_\l$, $\lar v_\l$, $v_\l\rar$, $v_\l\lar$.\par
\ni{\sl Definition 2.7.} For a module $V$, let $W_V$ be a set of
primitive vectors of $V$ corresponding to a composition series. We can
associate $W_V$ with a diagram defined by (2.1) for $v_\l,v_\mu\in W_V
$. We call this diagram the {\sl weight diagram} of $V$, and denote it
again by $W_V$.\qed \par From this definition, we see that the weight
diagram depends on the choices of primitive vectors: a module $V$ may
correspond to more than one weight diagrams. However, a weight diagram $
W_V$ does determine the structure of $V$ as we will see later. When
there is no confusion, we sometimes may use $V$ to mean its diagram or
vice versa. \par\ni{\sl Definition 2.8} A {\sl cyclic} module $X(\L)$
is a module generated by a primitive vector $v_\L$.\qed\par\ni{\bf
Theorem 2.9.} Suppose $V$ is an indecomposable module with a primitive
weight $\L$ and a primitive vector $v_\L$. We have\nl (1) If $\L$ is
typical, then $V=\VBL$.\nl (2) If $\L$ is multiply atypical, then all
primitive weights are multiply atypical.\nl (3) If $\L$ is singly
atypical, then all the primitive weights are singly atypical and $P_V=
\phi_\L^{(ij)}$ for some $i\le 0\le j$. \nl(4) For any choice of $W_V$,
$W_V$ must be connected, i.e., for any $u,v\in W_V$, there exist $u_0=u,
u_1,..,u_k=v\in W_V$ for some $k$ such that $u_i$ is linked to $u_{i+1}
$ by a line with an arrow for $i=0,...,k-1$.\par\ni{\sl Proof.} Take a
composition series: $0=V_0\subset V_1\subset...\subset V_n=V$. We will
prove the result by induction on $n$. If $n\le2$, by Lemma 2.4, we have
the result. Now suppose $n\ge 3$. We first prove {\sl Statement (A)}:
All primitive weights must have the same type: typical, singly atypical
or multiply atypical, and $W_V$ must be connected. {\sl Case (a):} $V/V_
1$ is decomposable. Decompose it into a direct sum of indecomposable
submodules: $V/V_1=\oplus_{i=1}^k V'_i/V_1$, then each $V'_i\supset V_1
$ must be indecomposable (if $V'_i=V''_i\oplus V_1$, we would then
write $V$ as a disjoint sum $V''_i\oplus(\Si_{j\ne i}V'_j)\,$), and by
inductive assumption, for all $i$, all primitive weights of $V'_i$,
which contains $V_1$, have the same type and $W_{V'_i}$ is connected,
and thus, all primitive weights of $V$ have the same type as that of $V_
1$ and $W_V$ is connected. {\sl Case (b):} $V/V_1$ is indecomposable.
By inductive assumption, all primitive weights in $V/V_1$ are of the
same type and $W_{V/V_1}$ is connected. This means that $(n-1)$ of $n$
primitive weights of $V$ are of the same type. If $V_1$ is the only
simple submodule of $V$, then $V_{n-1}$ is also indecomposable, and so
by inductive assumption, primitive weights of $V_{n-1}$ have the same
type and $W_{V_{n-1}}$ is connected; or else, if $V$ has another simple
submodule $V'_1$, then again by inductive assumption, primitive weights
of $V/V'_1$ have the same type and $W_{V/V'_1}$ is connected. In either
aspect, we can find another primitive vector of $V$ linking to, and
having the same type with, the primitive vector of $V_1$. Therefore all
primitive weights of $V$ have the same type and $W_V$ is connected.
This proves Statement (A), which implies (2)\&(4).\par Now (1) can be
proved as above by induction on $n$. To complete the proof of (3), by
Statement (A), we see that all primitive weights are now singly
atypical. We divide the proof into three steps: {\sl Step (i):} Suppose
$V$ is cyclic. Take $U_1=V/V_1$ and take $U_2=V_2$ if $V_1$ is the only
simple submodule of $V$, or else, $U_2=V/V'_1$ if $V'_1$ is another
simple submodule $V'_1$ of $V$. Then in either case, $U_1$ and $U_2$
are both cyclic. By inductive assumption, $P_{U_1}=\phi_\l^{(ij)}$, $P_
{U_2}=\phi_\mu^{(i'j')}$ for some $\l$, $\mu$, and $P_V$ is their
union. As they are not disjoint, and $\L$ is at least in one of them.
We see that $P_V$ has the required form. {\sl Step (ii):} For a
primitive weight $\l$, let $V^{(\l)}$ be the submodule generated by
primitive vectors with weights in $\phi_\l$. By Step (i), $P_{V^{(\l)}}
\subset\phi_\l$. For any two weights $\l,\mu$, \, $\phi_\l$ and $\phi_
\mu$ are either the same set or they are disjoint, thus different $V^{(
\l)}$ are disjoint. Since $V$, being a disjoint sum of $V^{(\l)}$, is
indecomposable, we must have $V=V^{(\L)}$. {\sl Step (iii):} Suppose $
\L^{(k)}\not\in P_V$ but $\L^{(i)},\L^{(j)}\in P_V$ for some $i<k<j$.
Let $W_1$ and $W_2$ be submodules generated by primitive vectors with
weight $\L^{(r)}$ such that $r>k$ and $r<k$ respectively. Then $V=W_1
\oplus W_2$, a contradiction with that $V$ is indecomposable. Thus we
have (3).\qed\par \ni{\sl Remark 2.10.} (1) If the weight diagram $W_V$
is connected, it is not necessary that $V$ is indecomposable (see
Remark 3.7). However, if $V$ is decomposable, we can always choose some
suitable primitive vectors such that $W_V$ is not connected.\nl (2) $v_
\mu\in U(G)v_\l$ does not mean that $v_\mu$ is linked to $v_\l$ as we
will see in Lemma 3.5(1.v).\qed\par\ \par\ \par\ni{\bf III.
INDECOMPOSABLE WEIGHT MODULES}\par Now we can classify indecomposable
modules $V$ with a primitive weight $\L$ of typical or singly atypical
type. We do this by classifying the weight diagram. If $\L$ is typical,
then by Theorem 2.9, $V$ is simply $\VBL$. Thus, from now on, we
suppose $\L$ is singly atypical. Then, all primitive weights are singly
atypical (again by Theorem 2.9).\par\ni{\sl Lemma 3.1.} For any
primitive vector $v_\l$, there exists at most one primitive vectors $v_
\mu$ such that $\mu<\l$ and $v_\mu\in U(G)v_\l$ (or $\mu>\l$ and $v_\mu
\in U(G)v_\l$).\par \ni{\sl Proof.} Suppose conversely there exists a
cyclic module $V_1=U(G)v_\l$ of the lowest dimension such that there
are 2 primitive vectors $v_\mu,v_\si \in V_1$ with, say, $\mu,\si<\l$.
By Theorem 2.9(3), $\mu,\si\le \l^-$. If there is a primitive vector $v_
\d\in V_1$ such that $\d>\l$ (and then $\d\ge\l^+$), then by our choice
of $V_1$ being lowest dimensional, $V'_1=U(G)v_\d$ does not have more
than one primitive weight $\d_1<\d$. Hence, by Theorem 2.9(3), if $\d_1<
\d$ is a primitive weight of $V'_1$, then $\d_1=\d^->\l^-\ge\mu,\si$.
Thus, if we let $V_2$ be the module generated by $\{v_\d\,|\,\d>\l,v_\d
\in V_1$ primitive$\}$, then $\mu,\si$ are not primitive weights in $V_
2$, so $v_\mu,v_\si$ are still primitive in $V_3=V_1/V_2$. Since $V_3$
is also a cyclic module, by our choice of $V_1$, we must have $V_2=0$,
i.e., $V_1$ is a highest weight module, but by Lemma 2.4(2.i), $V_1$
cannot contain 2 primitive weights $\mu,\si$. Thus we obtain a
contradiction.\qed\par\ni{\sl Corollary 3.2.} $W_V$ does not contain
(1) $u\rar v\rar w$, (2) $u^{\rar v}_{\rar w}$, (3) $u\lar v\lar w$, or
(4) ${}^{u\lar}_{v\lar}w$.\par\ni{\sl Proof.} This follows immediately
from Lemma 3.1 (say $V_W$ contains (1) or (2), then $v,w\in U(G)u$ and
their weights are less than the weight of $v$).\qed\par\ni{\sl
Corollary 3.3.} (1) If $v_\l\rar v_\mu$ or $v_\l\lar v_\mu$, then $\l=
\mu^+$. (2) $P_{U(G)v_\l}\subset\{\l^+,\l,\l^-\}$.\par\ni{\sl Proof.}
(1) Say, $v_\l\rar v_\mu$. Let $V_1=U(G)v_\l$. By Theorem 2.9(3), if $
\mu\ne\l^-$, then $\mu\le\l^{(-2)}$, and $\l^-$ must also be a
primitive weight of $V_1$. This contradicts Lemma 3.1. \ (2) follows
from (1) and Lemma 3.1.\qed\par\ni{\sl Lemma 3.4.} If (1) $v_{\L^+}\,^{
\lar\,v_\L}_{\rar\,u_\L}$ or (2) $^{v_\L\,\rar}_{u_\L\,\lar}\,v_{\L^-}$
is a part of the diagram, then $W_V$ must contain $X_4(\L)$: $v_{\L^+}
\,^{\lar\,v_\L\,\rar}_{\rar\,u_\L\,\lar}\,v_{\L^-}\,$. \par\ni{\sl
Proof.} Suppose conversely $V_1=U(G)v_\L$ is cyclic with the lowest
dimension such that its weight diagram contains, say, (2), but does not
contain $X_4(\L)$. By Lemma 2.4(2.i), $v_\L$ is not strongly primitive,
i.e., there exists $v_{\L^+}$ such that $v_{\L^+}\lar v_\L$. By
Corollary 3.2, we do not have a primitive vector $x$ such that $x\lar v_
{\L^+}\lar v_\L$, so $v_{\L^+}$ must be strongly primitive. Let $V_2=U(
G)v_{\L^+}\,$. If $u_\L\not\in V_2$, then the weight diagram of $V_1/V_
2$ contains $^{v_\L\,\rar}_{u_\L\,\lar}\,v_{\L^-}$ which violates the
choice of $V_1$ being of lowest dimension. Thus $u_\L\in V_2$, i.e., $u_
\L\in G_{-1}U(G_{-1})v_{\L^+}$, i.e., $v_{\L^+}\rar u_\L$, a
contradiction.\qed\par\ni{\sl Lemma 3.5.} If $V=X(\L)$ is cyclic,
then\nl(1) $W_V$ is one of (i) $v_\L$, (ii) $v_{\L^+}\lar v_\L$, (iii) $
v_\L\rar v_{\L^-}\,$, (iv) $v_{\L^+}\lar v_\L\rar v_{\L^-}$, or (v) $X_
4(\L)$. $X_4(\L)$ is the only cyclic module containing a primitive
vector which is not linked by the generator $v_\L$.\nl (2) $V$ is
uniquely (up to isomorphisms) determined by its diagram. \par\ni{\sl
Proof.} (1) If $v_\L$ is not linked to any primitive vector on one
side, then $V$ is a quotient of Kac- (or anti-Kac-) module, and we have
one of (i), (ii) or (iii). Suppose now $v_\L$ is linked to primitive
vectors on both sides, then $W_V$ contains (iv). If it is not (iv),
then there is another primitive vector $u$, by Lemma 3.1, it must have
weight $\L$ (and we denote $u_\L=u$). Since both modules $V/U(G)v_{\L^+}
$, $V/U(G)v_{\L^-}$ are Kac- or anti-Kac- modules, we must have $u_\L
\in U(G)v_{\L^+}$ and $u_\L\in U(G)v_{\L^-}$, i.e., $v_{\L^+}$ and $v_{
\L^-}$ are both linked to $u_\L$ with an arrow pointed to $u_\L$. By
Corollary 3.2, $u=u_\L$ is unique, i.e., we have (v). The statement
about $X_4(\L)$ is clearly true.\nl (2) Let $g_\L^+,g_\L^-$ be as in
Definition 2.6, if necessary, by replacing $v_{\L^+}$ by $g^+_\L v_\L$,
$v_{\L^-}$ by $g^-_\L v_\L$, we can suppose $v_{\L^+}=g^+_\L v_\L$, \ $
v_{\L^-}=g^-_\L v_\L$. This uniquely determines $V$ in the first 4
cases. Suppose now $V$ is (v), we can choose $u_\L=g^-_{\L^+}v_{\L^+}\,
$ and $g^+_{\L^-}v_{\L^-}$ must be a nonzero multiple of $u_\L$, i.e., $
g^+_{\L^-}v_{\L^-}=xu_\L$ for some $0\ne x\in\C$. Suppose $V'$ is
another module with weight diagram (v) (denote its corresponding
primitive vectors by the same notation with a prime) and suppose $g^+_{
\L^-}v'_{\L^-}=yu'_\L$ for some $0\ne y\in\C$, $y\ne x$. Form a direct
sum $V\oplus V'$ and let $V''$ be its submodule generated by $v''_\L=v_
\L+v'_\L$. Then $V''$ is indecomposable, which contains primitive
vectors $v''_{\L^+}=v_{\L^+}+v'_{\L^+}$, $v''_{\L^-}=v_{\L^-}+v'_{\L^-}
$, $u''_\L=u_\L+u'_\L$ and $u'''_\L=x u_\L+y u'_\L$. This contradicts
(1). (We will see from (4.11) that $x=-1$.)\qed\par In the following,
to be consistent, we use the same notations as in Ref. 4.\par\ni{\bf
Theorem 3.6.} (1) $W_V$ must be one of the following: (i) $X_4(\L)$,\nl
(ii)(a) $X_{5a}(\L,n):v_\L\rar v_1\lar v_2\rar...$ (ended by $v_{n-1}
\rar v_n$, or $v_{n-1}\lar v_n$ if $n$ is odd or even),\nl(b) $X_{5b}(
\L,n):v_\L\lar v_1\rar v_2\lar...$ (ended by $v_{n-1}\lar v_n$ or $v_{n-
1}\rar v_n$ if $n$ is odd or even),\nl where $v_i$ has weight $\L^{(-i)}
$.\nl (2) $V$ is uniquely determined up to isomorphisms by $W_V$.\par
\ni{\sl Proof.} For (2), Lemma 3.5(2) tells that $X_4(\L)$ determines $
V$. As in the proof of Lemma 3.5(2), $X_{5a}(\L,n)$, $X_{5b}(\L,n)$
also uniquely determine $V$. In fact, we can choose $v_i$ inductively,
such that $v_0=v_\L$, $v_{i-1}=g^+_{\L^{(-i)}}v_i$ if $v_{i-1}\lar v_i$
and $g^-_{\L^{(-i+1)}}v_{i-1}=v_i$ if $v_{i-1}\rar v_i$.\nl (1) ({\sl
cf.} proof of Theorem 4.9(2.ii).) If $V$ has a submodule corresponds to
(i), denote it by $V_1$; otherwise, let $V_1$ be a maximal submodule
whose diagram is (ii)(a) or (ii)(b). If $V$ has more primitive vectors,
choose $B$ to be its other primitive vectors such that $\Si_{v\in B}
$dim$(V_1\cap U(G)v)$ is minimum. If $B$ is not empty, by Theorem
2.9(4), there exists $v\in B$ linking to $V_1$. (i) If $V_1=X_4(\L)$.
By Corollary 3.2, the only possible links are: (a) $v_{\L^+}\lar v$,
(b) $v\rar v_{\L^-}$, (c) $v_{\L^+}\lar v\rar v_{\L^-}$, (d) $v\rar u_
\L$, (e) $u_\L\lar v$. For the first three cases, $v$ has weight $\L$.
Say, we have (c). We can choose $v$ such that $v_{\L^+}=g^+_\L v$ and $
g^-_\L v=x v_{\L^-}$, for some $x\ne0$. Apply Lemma 3.5(2) to $U(G)v$,
whose diagram is $X_4(\L)$, we must have $x=1$. Thus if we let $v'=v-v_
\L$, then $V_1\cap U(G)v'=0$. By replacing $v$ by $v'$ in $B$, we get a
contradiction with the choice of $B$. Similarly, for other cases, we
can also choose some primitive vector $v''\in V_1$ such that if we
replace $v$ by $v-v''$ we obtain a contradiction with the choice of $B
$. Therefore, $B$ is empty, and $V=X_4(\L)$. (ii) As $V_1$ is maximal,
by Lemma 3.4 and Corollary 3.2, $v$ can not be linked to $v_\L$ or $v_n
$. Also by Corollary 3.2, $v$ can not be linked to a vector $v_{2i}$ of
$X_{5a}(\L,n)$, or $v_{2i-1}$ of $X_{5b}(\L,n)$. On the other hand, if $
v$ is linked to a vector $v_{2i-1}\in X_{5a}(\L,n)$ or $v_{2i}\in X_{5b}
(\L,n)$, then $v$ must be linked to that vector with an arrow pointed
to it. Then as in (i), by replacing $v$ by $v-v''$ for some $v''\in V_1
$, we can get a contradiction. Thus again, $B$ is empty and we have
(ii).\qed\par \ni{\sl Remark 3.7.} Diagrams such as $u^{\lar v}_{\lar w}
$ and $^{v\rar}_{w\rar}u$ can exist, but they correspond to
decomposable modules: by replacing $w$ by $w-v$, we see that $w$ is not
linked to $u$, $v$ (see Remark 2.10). \qed\par It is not difficult to
construct $X_{5a}(\L,n)$, $X_{5b}(\L,n)$ as follows: in the module ${
\ol V}_*((\L^-)_{low}))\oplus{\ol V}(\L^-)$, whose diagram has two
parts: $\lar v'_{\L^-}$, $v''_{\L^-}\rar$, ``joining'' two primitive
vectors $v'_{\L^-}$, $v''_{\L^-}$ into one, by letting $v_{\L^-}=v'_{\L^
-}+v''_{\L^-}$, we obtain $X_{5b}(\L,3)=U(G)v_{\L^-}\,$. Now $X_{5b}(\L^
{(-2i)},3)$ has diagram $v'_{\L^{(-2i)}}\lar v_{\L^{(-2i-1)}}\rar v''_{
\L^{(-2i-2)}}\,$. By taking a quotient module, by ``merging'' $v'_{\L^{(
-2i)}}$, $v''_{\L^{(-2i)}}$ into one, we obtain $X_{5b}(\L,2n)=\oplus_
{i=0}^n X_{5b}(\L^{(-2i)},3)/\oplus_{i=1}^n U(G)(v'_{\L^{(-2i)}}-v''_{
\L^{(-2i)}})\,$. Modules $X_{5a}(\L,n),X_{5b}(\L,2n+1)$ can be realized
as subquotients of $X_{5b}(\L,k)$ for some $k$. To construct $X_4(\L)$,
form an induced module ${\wh V}(\L)=\mbox{Ind}_{G_0}^G V^0(\L)=U(G)
\otimes_{U(G_0)}V^0(\L)\cong U(G_{-1})\otimes U(G_{+1})\otimes V^0(\L)
$. (Note that this module is in general decomposable, therefore not
cyclic. However a cyclic module can be realized as its quotient
module.) We are not going to realize $X_4(\L)$ to be a quotient module
of ${\wh V}(\L)$, but as a submodule of ${\wh V}(\L^+-2\rho_1)$: let $v_
{\L^+-2\rho_1}$ be the highest weight vector in $G_0$-module $V^0(\L^+-
2\rho_1)$, then $v_{\L^+}=gv_{\L^+-2\rho_1}$ (where $g$ is the highest
root vector in $U(G_{+1})$\,) is a strongly primitive vector with
weight $\L^+$, which generates Kac-module ${\ol V}(\L^+)$. It is clear
that there must exist a primitive vector $v_\L$ in ${\wh V}(\L^+-2\rho_
1)$ such that $v_{\L^+}\lar v_\L$. By Lemma 3.4 and Theorem 3.6, $v_\L$
generates a module corresponding to $X_4(\L)$ (see also (4.11)).\par We
see that just as an anti-Kac-module is isomorphic to a Kac-module, $X_
{5b}(\L,2n+1)$ is isomorphic to $X_{5a}(\Si,2n+1)$ for some $\Si$ (in
fact, $(\Si^{(-2n-1)})_{low}=-\L$). To see that $X_4(\L)$, $X_{5a}(\L,n)
$, $X_{5b}(\L,n)$ are indecomposable: suppose $V=V_1 \oplus V_2$ is a
disjoint sum, then each simple submodule must be contained in $V_1$ or $
V_2$, and then we can obtain that all primitive vectors must be in one,
say, $V_1$, and $V=V_1$. Now we can conclude the following\par\ni{\bf
Theorem 3.8.} \, $\{\VBL\,|\,\L$ typical$\}\cup\{X_4(\L),X_{5a}(\L,n),X_
{5b}(\L,2n)\,|\,n\in\Z_+\bs\{0\}\,\}$ \,is \,the complete set of
indecomposable modules with at least a primitive weight of typical or
singly atypical type.\par\ni{\sl Proof.} It remains to prove there is
no isomorphism between each other. This can be seen by comparing number
of simple submodules and number of composition factors.\qed\par\ \par\
\par\ni{\bf IV. INDECOMPOSABLE GENERALIZED WEIGHT MODULES}\par Now
suppose $V$ is an indecomposable $G$-module such that $H$ acts on $V$
not necessarily diagonally. Such a module is called a generalized
weight module (in Ref. 4, it is called a module with nondiagonal Cartan
subalgebra, or nondiagonal module). In this case, we do not have weight
space decomposition. However, by the properties of semi-simple Lie
algebras, we see that $V$ must be $H_0$-diagonal, where $H_0$ is the
Cartan subalgebra of $G'_0$ (where, here and after, $G'_0$ is the
subalgebra of $G_0$ with co-dimension 1 such that $h_0\not\in G'_0$)
with a basis $\{h_i\,|\,i\ne0\}$. We have $$V=\oplus_{\l\in H^*}{\bf V}_
\l,\ \mbox{ where }\ {\bf V}_\l=\{v\in V\,|\,(h-\l(h))^n v=0\mbox{ for
some }n\in\Z_+\}.\eqno(4.1)$$ ${\bf V}_\l$ are called {\sl generalized}
weight spaces. With this decomposition, we have, similar to Sect. II,
notions of {\sl generalized} (strongly) (anti-) primitive vectors
(weights), {\sl generalized} weight diagram, etc. In the following, we
often omit the word {\sl generalized} if there is no confusion.\par We
can take a composition series $$0=V^{(0)}\subset V^{(1)}\subset...
\subset V^{(n)}=V,\eqno(4.2)$$ such that each $V^{(i)}$ is a direct sum
of subspaces $$V^{(i)}=V^{(i-1)}\oplus{\ol V}^{(i)},\ \mbox{ where }\ {
\ol V}^{(i)}\mbox{ is a $G'_0$-module.}\eqno(4.3)$$ Then each ${\ol V}^
{(i)}$ has a unique, up to scalars, generalized primitive vector $v^{(i)
}_{\l_i}$ with weight $\l_i$. Sometimes, we can just choose (4.2) to be
any series of submodules such that (4.3) holds.\par\ni{\sl Definition
4.1.} For an integral dominant weight $\L$, construct an indecomposable
module, the {\sl generalized Kac-module} ${\ol V}(\L,n)$, $n\in\Z_+\bs
\{0\}$, as follows: It is a semi-direct sum of $n$ copies of $\VBL$,
such that each copy is a $G'_0\oplus G_{-1}$-module, and $$h_0v^{(i)}_
\L=\L(h_0)v^{(i)}_\L+v^{(i-1)}_\L,i=2,...,n,\eqno(4.4)$$ where $v^{(i)}_
\L$ belongs to the $i$-th copy of $\VBL$. As $\VBL$ is the induced
module, it is easy to check that ${\ol V}(\L,n)$ is well defined as an
indecomposable $G$-module. Similarly, we can define generalized
anti-Kac-module ${\ol V}_*(-\L,n)$.\qed\par Note that if $\L$ is singly
atypical, then ${\ol V}(\L,n)$ has the ``zigzag'' weight diagram such
that $$v^{(i)}_\L\rar v^{(i)}_{\L^-},\ \ \ v^{(i-1)}_\L\lar v^{(i)}_{\L^
-},\eqno(4.5)$$for $i=n,n-1,...,1$ (where we take $v^{(0)}_\L=0$).\par
\ni{\sl Remark 4.2.} If $\L$ is typical, ${\ol V}(\L,n)$ give us
examples that there may exist a primitive vector $v_\L$ such that it is
not linked to any primitive vector, but $U(G)v_\L$ is not simple and
that an indecomposable module may not correspond to a connected weight
diagram. Therefore, more care should be taken when we consider
generalized weight modules. We shall see in Lemmas 4.5(2)\&4.6(3), that
this does not happen if $\L$ is singly atypical.\qed\par \ni{\sl
Definition 4.3.} Define a partial order on generalized primitive
vectors: we say $v$ has {\sl higher level} than $u$ or $v$ is on {\sl
top} of $u$ $\ \Rla\ $ $u\in U(G)v$, but $v\not\in U(G)u$.\qed\par We
see that the generalized primitive vectors of ${\ol V}_0(\L,n)={\ol V}(
\L,n)$ are well ordered by this partial order. Now by removing top
vector $v^{(n)}_\L$, removing bottom vector $v^{(1)}_{\L^-}$, and
removing both top and bottom vectors $v^{(n)}_\L$, $v^{(1)}_{\L^-}$
respectively, we obtain three indecomposable modules: ${\ol V}_1(\L,n)
${}$=${}$U(G)v^{(n)}_{\L^-}$, ${\ol V}_2(\L,n)${}$=${}${\ol V}(\L,n)/U(
G)v^{(1)}_{\L^-}$ and ${\ol V}_3(\L,n-1)${}$=${}$U(G)v^{(n)}_{\L^-}/U(G)
v^{(1)}_{\L^-}\,$. We see that generalized anti-Kac-module ${\ol V}_*((
\L^-)_{low},n)$ can be realized as ${\ol V}_3(\L,n)$, a subquotient of
generalized Kac-module ${\ol V}(\L,n+1)$.\par\ni{\sl Lemma 4.4.}
Suppose $V$ is a (generalized) highest weight module with highest
weight $\L$ (i.e., $V$ is generated by a generalized strongly primitive
vector $v_\L$). (1) If $\L$ is typical, then $V\cong {\ol V}(\L,n)$ for
some $n$. (2) If $\L$ is singly atypical, then $V$ is a quotient of ${
\ol V}(\L,n)$ for some $n$. More precisely, $V={\ol V}_i(\L,n)$, $i=0$
or 2.\par\ni{\sl Proof.} First, as $U(G)=U(G_-)U(H)U(G_+)$, we see that
$V=U(G_-)U(H)v_\L$ does not have weight $>\L$. Let $V$ be as in (4.2).
We use induction on $n$. If $n=1$, the result is obvious. Suppose now $
n\ge2$. As $V'=V/V^{(1)}$ is still a highest weight module, by
inductive assumption, it has the required form. If $\L$ is typical,
then $V'={\ol V}(\L,n-1)$, and therefore by (4.4), we can inductively
choose spaces ${\ol V}^{(i)}$ in (4.3) as a copy of ${\ol V}(\L)$ with
the primitive vector $v^{(i)}_\L$, $i=n,n-1,...,2$, such that $$h_0 v^{(
i)}_\L=\L(h_0)v^{(i)}_\L+v^{(i-1)}_\L+v_i\ \mbox{ for some }\ v_i\in V^
{(1)},\eqno(4.6)$$ where, we take $v^{(1)}_\L$=0 when $i=2$. As all $v^
{(i)}_\L\in {\bf V}_\L$, we have $v_i\in{\bf V}_\L$, i.e., $v_i$ has
weight $\L$. By replacing $v^{(i-1)}_\L$ by $v^{(i-1)}_\L+v_i$ (and
replacing space ${\ol V}^{(i)}$ by $U(G_{-1})U(G'_0)(v^{(i-1)}_\L+v_i)$
accordingly), we can suppose $v_i=0$ if $i>2$. If $v_2\ne0$, denote it
by $v^{(1)}_\L$, then we have (4.4), and thus $V={\ol V}(\L,n)$ and the
result follows. On the other hand, if $v_2=0$, then $V''=\Si_{i=2}^n {
\ol V}^{(i)}$ is a submodule of $V$ and $V=V^{(1)}\oplus V''$, a
contradiction with that $V$ is a highest weight module. This proves
(1). For (2), as $V'$ has the required form, without loss of
generality, say $V'={\ol V}_2(\L,n)$. Then as space, it is the direct
sum of $(n-1)$ copies of $\VBL$ plus $V(\L)$. Now again choose ${\ol V}^
{(i)}$ to be a copy of $\VBL$ for $i\ge3$ and ${\ol V}^{(2)}=V(\L)$ (in
this case, (4.2) is not a composition series). Now follow the arguments
exactly as above, we have (4.4) for $i>=3$. Take $v^{(i)}_{\L^-}$ to be
the other primitive vector in ${\ol V}^{(i)}$ for $i\ge 3$. Then we
have (4.5) for $i=n,..,3$. Now we must have ${\ol V}^{(1)}\subset V'''=
U(G)v^{(2)}_\L$ (otherwise, $V={\ol V}^{(1)}\oplus(\Si_{i=2}^n{\ol V}^{(
i)})$ is decomposable). It remains to prove $V'''$, which is now as
space ${\ol V}^{(1)}\oplus {\ol V}^{(2)}$, is $\VBL$ (and then $V={\ol
V}(\L,n)$\,). This follows from Lemma 4.5(3) below.\qed\par The
following Lemma 4.5 tells that, unlike typical case, for atypical
weight $\L$, we do not have indecomposable module whose composition
factors are $n$ copies of $V(\L)$.\par\ni{\sl Lemma 4.5.} (1) Suppose $
V$ is a (generalized weight) module whose composition factors are $n$
copies of $V(\L)$ with $\L$ typical or singly atypical, then (i) if $V$
is indecomposable, then either $n=1$, or else, $\L$ is typical and $V={
\ol V}(\L,n)$; (ii) if $\L$ is singly atypical, then $V$ is the direct
sum of $n$ copies of $V(\L)$ (and thus, $V$ is a weight module).\nl(2)
If $v_\L$ is a primitive vector with singly atypical weight $\L$ such
that it is not linked to any primitive vector, then $U(G)v_\L$ is
simple.\nl(3) If $V$ is an indecomposable module with 2 composition
factors of singly atypical type, then $V$ is a Kac- (or anti-Kac-)
weight module.\par\ni{\sl Proof.} (1) Let (4.2) be a composition series
of $V$. By induction on $n$, we see that we only need to prove the
result for $n=2$. (i) If $\L$ is typical, it is easy to see, as in the
proof of Lemma 4.4, $V={\ol V}(\L,2)$. (ii) Suppose now $\L$ is singly
atypical, let $v^{(i)}_\L$, $i=1,2$, be the primitive vectors of $V$
and suppose $$h_0 v^{(1)}=\L(h_0)v^{(1)},\ h_0 v^{(2)}_\L=\L(h_0)v^{(2)}
_\L+av^{(1)}_\L,\ \mbox{ for some }\ a\in\C.\eqno(4.7)$$ Note that, in $
U(G)$, for $G=sl(m/n)$, using notations in Sect. II, we have $$\prod_{-
m\le i\le0\le j\le n}e_{ij}\ \prod_{-m\le i\le0\le j\le n}f_{ij}=\si
\prod_{-m\le i\le0\le j\le n}(\Si_{k=i}^0 h_k-\Si_{k=1}^j h_k-i-j)+g^+,
\eqno(4.8)$$ for some $g^+\in U(G)G^+$, where $\si=\pm1$ (this can be
proved by ordering $e_{ij}$ ($f_{ij}$) properly in the products, such
that if $j-i>j'-i'$, or $j-i=j'-i'$ and $j>j'$, then $e_{ij}$ ($f_{ij}
$) is placed to the right (left) of $e_{i'j'}$ ($f_{i'j'}$ resp.); then
using induction on $m,n$). We have $$U(G)G^+v^{(2)}_\L=0,\ \ \ \ \
\prod_{-m\le i\le0\le j\le n}f_{ij}v^{(2)}_\L=0.\eqno(4.9)$$ (The {\sl
l.h.s.} of the second equality is in the bottom composition factor $V(
\L^-)$ of $\VBL$ with weight $\L-2\rho_1$, any copy of $V(\L)$ does not
have a vector with weight $\L-2\rho_1$; therefore it is zero in $V$).
Now apply (4.8) to $v^{(2)}_\L$, using (4.9) and (4.7), we obtain $$0=
\prod_{-m\le i\le0\le n}(\L+\rho,\a_{ij})v^{(2)}_\L+\ \sum_{-m\le i\le0
\le j\le n}\ \ \prod_{-m\le i'\le0\le j'\le n,\,(i'j')\ne(i,j)}(\L+\rho,
\a_{i'j'})av^{(1)}_\L.\eqno(4.10)$$ However, by Definition 2.2, there
is exactly one atypical root, i.e., one pair of $(i,j)$ such that $(\L+
\rho,\a_{ij})=0$, thus (4.10) forces $a=0$. This proves $V=U(G)v^{(2)}_
\L\oplus U(G)v^{(1)}_\L=V(\L)\oplus V(\L)$.\nl(2) Let $V=U(G)v_\L$ be
as in (4.2). If $n>1$, using induction, we can suppose $V/V^{(1)}$ is
simple, i.e., $n=2$. By (1), two composition factors can not be the
same, but then as in (4.6) (for $i=2$) and the arguments after (4.6),
we see that $U(G)v_\L$ is a weight module. As $v_\L$ is not linked to
any primitive vector, by Theorem 2.9, $V$ has to be simple (and so $n=2
$ does not occur).\nl(3) The proof is the same as (2).\qed\par Now we
have the generalized weight modules ${\ol V}_i(\L,n)$ and $X_{5a}(\L,n)
$, $X_{5b}(\L,n)$, whose diagrams have 2 {\sl endpoints} (i.e., vectors
linked by only one vector). Such diagrams are called {\sl lines}. We
can ``join'' and ``merge'' those modules to form other indecomposable
modules just as we did to form $X_{5b}(\L,n)$ in Sect. III. In
particular, we can define a module $X_4(\L,m,n,x)$ ($m,n\ge 2$, $0\ne x
\in\C$) as follows: let $X_3(\L,m,n)=({\ol V}_1 (\L^+,m)\oplus{\ol V}_2(
\L,n))/U(G)(v^{(1)}_\L-v'^{(1)}_\L)$ (with obvious meanings of
notations) be the quotient module (by ``merging'' the two bottom
endpoints $v^{(1)}_\L\in{\ol V}_1(\L^+,m)$ and $v'^{(1)}_\L\in{\ol V}_2(
\L,n)$\,). Then $X_4(\L,m,n,x)$ is the submodule of $X_3(\L,n)$
generated by $v^{(m)}_\L+xv'^{(n)}_\L$ (by ``joining'' the two top
endpoints). By this construction, one sees that $X_4(\L,m,n,x)$ is
indecomposable, generated by a primitive vector, therefore cyclic. Its
diagram is a {\sl circle}, i.e., no endpoints. It is interesting to see
that $X_4(\L)$ can be realized as $$X_4(\L)=X_4(\L,2,2,-1).\eqno(4.11)
$$ This is because: if $x=-1$, when we ``join'' the top endpoints, the
second term of the {\sl r.h.s.} of (4.4) is lost, and so $h_0$ becomes
acting diagonally. We point out that only with a circle, an $x$ makes
difference: just as in Sect. III, we can choose suitable primitive
vectors starting from a vector $v_{\l_1}$ of the weight diagram (we
always choose $v_{\l_1}$ to be an endpoint if it is a line), and follow
the links between vectors, such that if $v_\l\rar v_\mu$, then $g^-_\l
v_\l=v_\mu$ and if $v_\l\lar v_\mu$, then $v_\l=g^+_\mu v_\mu$. But
within a circle, the last vector $v_{\l_n}$ we chose is linked to the
first one and in this case, say, $v_{\l_n}\lar v_{\l_1}$, we may have $
g^+_{\l_1}v_{\l_1}=xv_{\l_n}$ for some $x\in\C\bs\{0\}$. By rescaling
vectors, we see that $x$ can be shifted anywhere, but can not be
eliminated if the diagram is a circle.\par It is interesting to see
that we can ``add'' a primitive vector to a circle to break it into a
line: in the above construction, $X_3(\L,m,n)$ can be obtained by
adding $v^{(m)}_\L$ to $X_4(\L,m,n,x)$.\par\ni{\sl Lemma 4.6.} (1) Let $
V$ be a cyclic module with a primitive weight $\L$ of typical or singly
atypical type, then all generalized primitive weights have the same
type. Furthermore, $V$ must be (i) ${\ol V}(\L,n)$ if $\L$ is typical
or (ii) a quotient of some $X_4(\Si,m,n,x)$, where $\Si=\L^{(i)}$ for
some $i=-1,0,1$ (and so its diagram is obtained from $X_4(\Si,m,n,x)$
by removing some vectors from the bottom) if $\L$ is singly atypical.\nl
(2) If $\L$ is singly atypical and $V=X(\L)$ is cyclic, then $P_V\in\{
\L^+,\L,\L^-\}$.\nl(3) If $V$ is an indecomposable module with a singly
atypical weight, then $W_V$ is connected.\par\ni{\sl Proof.} (1) Take a
composition series $0=V_0\subset V_1\subset...\subset V_k=V$. First as
in the proof of Theorem 2.9, by induction on $k$, we can prove all
generalized primitive weights have the same type. Now if $k=1$, we
clearly have the result. Suppose $k\ge 2$. Then $V/V_1$ is still
cyclic, by inductive assumption, it has the required form. Now follows
exactly the same arguments as in the proof of Lemma 4.4, we have the
result. (2) follows from (1). \ (3) By (1), diagrams for cyclic modules
are connected. If $W_V$ is not connected, two submodules generated by
primitive vectors in two disjoint parts are disjoint and $V$ is their
disjoint sum.\qed\par In ${\ol V}(\L,n)$, all $v^{(i)}_{\L^-}$ have
weight $\L^-<\L$, thus, Lemma 3.1 does not hold for generalized weight
modules; and we have $v^{(i)}_{\L^-}\in G_{-1}U(G^-)U(H)v_\L^{(n)}$,
but $v^{(i)}_{\L^-}\not\in G_{-1}U(G^-)v_\L^{(n)}$ if $i<n$, thus by
(2.1), $v^{(i)}_{\L^-}$ is not linked by the generator $v_\L^{(n)}$ if $
i<n$. However, from the structure of weight diagrams of cyclic modules
in Lemma 3.6, we see that Corollary 3.2 still holds for generalized
weight diagrams. Now similar (but not exactly) to Ref. 4, we define\par
\ni{\sl Definition 4.7.} (1) A {\sl chain} is a triple $(C,\L,x)$,
where $\L$ is a singly atypical weight, $C$ is a finite set: $C=\{v_1,..
 .,v_n\}$, $x\in\C$, such that $v_i$ has a weight $\l_i=\L^{(k)}$ for
some $k$ and there are links with arrows between elements in $C$, such
that\nl(i) each element $v$ is linked by at most 2 elements, and in
this case, it can only be one of (a) $\lar v\rar$, (b) $\rar v\lar$,
(c) $v^\lar_\rar$, (d) ${}^\rar_\lar v$ (accordingly, $v$ is called a
{\sl top}, {\sl bottom}, {\sl left}, {\sl right}, point);\nl(ii) If $v_
i$ is linked to $v_j$ with an arrow pointed to it, then $v_i$ is not
derived from $v_j$ (we say $v$ is {\sl derived} from $u$ if there exist
$v=u_0,u_1,...,u_k=u$ such that $u_{i+1}$ is linked to $u_i$ with an
arrow pointed to it for $i=0,..,k-1$);\nl(iii) $C$ is connected;\nl
(iv) if $v_i$ is linked to $v_j$ from the left, then $\l_i=\l_j^+$;\nl
(v) if $v_i$ is a {\sl leftmost} element, i.e., there is no elements
linked to it from left side, then $\l_i=\L$;\nl(vi) $x=0\ \ \Rla\ \ $
there are 2 endpoints.\nl(2) A {\sl subchain} of $(C,\L,x)$ is a subset
of $C$ together with the weight, the original relationship between the
elements.\nl(3) An {\sl isomorphism} between two chains $(C,\L,x)$ and $
(C',\L',x')$ is a bijection $C\ \rar\ C'$ which preserves weight,
linking relationship and $x=x'$.\nl(4) An {\sl anti-isomorphism}
between $(C,\L,x)$ and $(C',\L',x')$ is a bijection $\phi:C\ \rar\ C'$
such that (i) $u\lar v\ \Rla\ \phi(v)\rar\phi(u)$, and $u\rar v\ \Rla\
\phi(v)\lar\phi(u)$; (ii) if a {\sl rightmost} element of $C'$ has the
weight $\l$, then $\l_{low}=-\L$; (iii) $x=x'$.\qed\par\ni{\sl Remark
4.8.} (1) For a chain $(C,\L,x)$, we can break it, at top and bottom
points, into pieces of subchains according to the rule: if $\lar v_\l
\rar$ (or $\rar u_\l\lar$), then we break at the point $v_\l$ (or $u_\l
$) into $\lar v'_\l,v''_\l\rar$ (or $\rar u'_\l,u''_\l\lar$, resp.).
Then each piece is corresponding to some ${\ol V}_i(\L^{(j)},k)$. If we
only break it at the bottom points, then each piece corresponds to a
cyclic module. This gives us a better understanding of what a chain
is.\nl(2) If a chain is a line, then conditions (ii)\&(iv) are
unnecessary as (ii) can not happen and we have a unique way to
associate with each $v_i$ a weight $\l_i$ such that (iv) is satisfied.
\nl(3) Examples of non-chains: (a) a diagram such that ${}^{v_1\,\rar}_
{v_3\,\lar}\,v_2$ and ${}^{v_1\,\lar}_{v_3\,\rar}\,v_4$ (it violates
(ii)); (b) $v_1\,^{\ssc\rar\lar}_{\dis\rar}{}^{\ssc\rar\lar}_{\dis\lar}
\,v_5$ (it violates (iv)).\qed\par We see that the weight diagrams of
all indecomposable modules introduced up to now are chains. Now we can
prove the main result of this section.\par\ni{\bf Theorem 4.9.} (1) For
any chain $(C,\L,x)$, there is a unique indecomposable generalized
weight $G$-module $X(C,\L,x)$ corresponding to it.\nl (2) If $V$ is an
indecomposable generalized weight module with a primitive weight $\L$
of typical or singly atypical type, then (i) $V={\ol V}(\L,n)$ if $\L$
is typical, or (ii) there exists a unique chain $(C,\L,x)$, up to an
(anti-) isomorphism, such that $V=X(C,\L,x)$ if $\L$ is singly
atypical.\par\ni{\sl Proof.} (1) Break the chain as in Remark 4.8, let $
V'$ be the direct sum of all ${\ol V}_i(\L^{(j)},k)$ obtained. Then for
each pair of $u'_\l,u''_\l$, we can ``merge'' them into one $u_\l$
(this is a bottom point) by taking quotient $V'/U(G)(u'_\l-u''_\l)$ and
let $V''$ be this quotient module; and for each pair of $v'_\l,v''_\l$,
we ``join'' them into one $v_\l$ (this is a top point) by letting $v_\l=
v'_\l+av''_\l$ (where $a=x$ if $x\ne0$ and $v_\l$ is the last vector to
be joined in order to form the circle; or otherwise, $a=1$). Now let $X(
C,\L,x)$ be the submodule of $V''$ generated by all ``joined'' vectors $
v_\l$ (all top points). Then we see $X(C,\L,x)$ is corresponding to $(C,
\L,x)$. By the statements following (4.11), we see that it is uniquely
determined by $(C,\L,x)$. To see it is indecomposable, suppose it is a
direct sum of submodules $V_1\oplus V_2$. If $v_\l$ plus a linear
combination of other vectors in $C$ with weight $\l$ belongs to $V_1$,
then by Lemma 4.6(1.ii), all vectors derived from $v_\l$ is in $V_1$;
as $C$ is connected, all vectors in $C$ must also be in $V_1$; thus $V_
2=0$.\nl(2)(i) Let $V'$ be the submodule generated by $\{v_\l\in C\,|\,
\l$ singly atypical (if any)$\}$, and let $V^{(\l)}$ be the submodule
generated by $\{v\in C\,|\,v$ has weight $\l\}$ if $\l$ is typical.
Then by Lemma 4.6(1.i), $V=\oplus_{\{\l\ {\rm typical}\}}V^{(\l)}\oplus
V'$ is a direct sum, thus $V=V^{(\L)}$. Now by Lemmas 4.6(1.i)\&
4.5(1.i), we have $V={\ol V}(\L,n)$.\nl(ii) ({\sl cf.} proof of Theorem
3.6.(1)) We want to prove {\sl Statement (A):} For any generalized
weight module $V$, we can choose $W_V$ such that it is a union of
pieces of disconnected subdiagrams, each subdiagram is a chain. Then by
Lemma 4.6(3), we have the result. Now let $V$ has a composition series
as in (4.2). If $n\le2$, (A) follows immediately from Lemma 4.5. Assume
now $n\ge3$. By induction, suppose $W_{V^{(n-1)}}$ satisfies (A). Let $
v$ be a primitive vector corresponding to $V/V^{(n-1)}$. Then $v$ is
not derived by any primitive vector. To understand it better, we prove
(A) case by case. Using Corollary 3.2, we consider all possible cases
below. {\sl Case (a):} If $v$ is not linked to $W_{V^{(n-1)}}$, then $v
$ itself is a piece of chain. Thus we have (A). {\sl Case (b):} If $v$
is linked to $W_{V^{(n-1)}}$ only on the right side, then there is a
unique $u$ in a piece of $W_{V^{(n-1)}}$ such that $u\lar v$. Thus,
there does not exist $u'$ such that $u'\lar u$. {\sl (b.1):} If there
is $v'\in W_{V^{(n-1)}}$ such that $u\lar v'$, then we can replace $v$
by $v-v'$ so that $v$ is now not linked to $W_{V^{(n-1)}}$. This
becomes (a). {\sl (b.2):} It remains that $u$ is either a left top
endpoint: $u\rar w$, or a right bottom endpoint: $w\rar u$. Then we see
that $v$ can be added to that piece so that we have $u^{\lar v}_{\rar w}
$ or $w\rar u\lar v$, and it is still a chain. Thus we have (A). {\sl
Case (c):} Similar to (b), if $v$ is linked to $W_{V^{(n-1)}}$ only on
the left side, we have (A). {\sl Case (d):} If $v$ is linked to $W_{V^{(
n-1)}}$ on both sides, i.e., there are unique $u,w\in W_{V^{(n-1)}}$
such that we have {\sl Diagram (D.1):} $u\lar v\rar w$. {\sl (d.1):} If
$u$, $w$ belong to 2 different pieces of chains, then as in (a) and
(b), by a suitable choice of $v$ and rescaling primitive vectors if
necessary, either $v$ can be added to one piece, or we can link 2
pieces through $v$ into one piece of chain. Thus we have (A). {\sl
(d.2):} Finally suppose $u$ and $w$ are in the same piece of chain. In
this case, suppose $v$ has weight $\l$, and $u=g^+_\l v$, $g^-_\l v=xw$
for some $x\ne0$. {\sl (d.2.1):} If there is $v'\in W_{V^{(n-1)}}$ such
that we have {\sl Diagram (D.2):} $u\lar v'\rar w$. {\sl (d.2.1.i):} If
$x=1$. Replace $v$ by $v-v'$, then $v$ is not linked to any one. This
becomes (a) and we have (A). {\sl (d.2.1.ii):} If $x\ne1$, by replacing
$v$ by $v-v'$ and replacing $v'$ by $xv-v'$, when we add $v$ into that
piece, from (D.1) and (D.2), we see that (D.2) is broken into: $u\lar
v'$, $v\rar w$ ({\sl cf.} the statement before Lemma 4.6). Therefore it
is still corresponding to a piece of chain (or 2 pieces of chains if
that piece becomes disconnected). Thus we have (A). {\sl (d.2.2):} If
there exists $v'\in W_{V^{(n-1)}}$ such that $v'\rar w$, but not $u\lar
v'$ or there exists $v''\in W_{V^{(n-1)}}$ such that $u\lar v''$ but
not $v''\rar w$, or if both, by replacing $v$ by $v-v'$, or $v-v''$, or
$v-v'-v''$, we see that (D.1) becomes $u\lar v$ or $v\rar w$, or $v$.
This becomes (b) or (c) or (a). {\sl (d.2.3):} It remains that $u$ is
either a right bottom endpoint $t\rar u$ or a left top endpoint $u\rar
t$ and $w$ is either a left bottom endpoint $w\lar x$ or a right top
endpoint $x\lar w$. From (D.1), we see that $v$ can be added to it so
that we have $t\rar u\lar v\rar w\lar x$ or $t\rar u^{\lar}{}^{v\rar}_
{x\lar}w$, or $u_{\rar t}^{\lar v\rar}w\lar x$ or $u_{\rar t}^{\lar}{}^
{\!v\,\rar}_{\,x\lar}w$, and it is still a chain. Thus we have (A).
This completes the proof of Statement (A).\qed\par\ni{\sl Remark 4.10.}
By applying the above results to $sl(m/1)$, as all primitive weights of
$sl(m/1)$ are typical or singly atypical, we have efficiently
classified all finite dimensional (weight or generalized weight) $sl(m/
1)$-modules.\qed\par\ \par\ \par\ni{\bf REFERENCE} \par\ni\hi2ex\ha1 $^
{1}$ V.G. Kac, ``Representations of classical Lie superalgebras,'' in
{\sl Lecture in Mathematics} {\bf 676}, Springer-Verlag (1977), ed.
Bleuler K., Petry H. and Reetz A., 579-626.\par\ni\hi2ex\ha1 $^{2}$ J.
Van der Jeugt, J.W.B. Hughes, R.C. King \& J. Thierry-Mieg, ``A
character formula for singly atypical modules of the Lie superalgebra $
sl(m/n)$'', {\sl Comm. Alg.} {\bf 18}(1990), 3453-3480. \par\ni\hi2ex\ha
1 $^{3}$ J.W.B. Hughes, R.C. King \& J. Van der Jeugt, ``On the
composition factors of Kac-modules for the Lie superalgebra $sl(m/n)
$,'' {\sl J. Math. Phys.} {\bf 33}(1992), 470-492. \par\ni\hi2ex\ha1 $^
{4}$ Y.-C. Su, ``Classification of finite dimensional modules of the
Lie superalgebra $sl(2/1)$,'' {\sl Comm. Alg.}, {\bf 20}(1992),
3259-3277.\par\ni\hi2ex\ha1 $^{5}$ V.G. Kac, {\sl Infinite-dimensional
Lie algebras}, Birkh\"auser Boston, Cambridge, MA, 2nd ed., \nopagebreak
1985.\end{document}